\documentclass[letter,10pt,conference]{style/ieeeconf}
\IEEEoverridecommandlockouts
\usepackage{amsmath,amsfonts,amssymb}
\usepackage{graphicx}
\usepackage{gensymb}
\usepackage{amsthm}
\usepackage{algpseudocode}
\usepackage[ruled,vlined,linesnumbered]{algorithm2e}
\newcommand\figcapskip{\vspace*{-.75em}}
\newcommand\figbottomskip{\vskip -1.5em}
\usepackage[utf8]{inputenc}
\usepackage{xcolor}
\newtheorem{theorem}{Theorem}
\newtheorem{lemma}[theorem]{Lemma}
\newtheorem{definition}[theorem]{Definition}

\newtheorem{assumption}[theorem]{Assumption}

\begin{document}

\title{\LARGE \bf
Near-optimal control of nonlinear systems with hybrid inputs and dwell-time constraints
\vspace*{-10pt}
}

\author{
Ioana Lal, Constantin Mor\u{a}rescu, Jamal Daafouz, Lucian Bu\c{s}oniu
\thanks{I.~Lal and L.~Bu\c soniu are with the Automation Department, Technical University of Cluj-Napoca, Romania. C. Mor\u{a}rescu and J. Daafouz are with Universit\'e de Lorraine, CRAN, UMR 7039 and CNRS, CRAN, UMR 7039, Nancy, France. Email addresses: \texttt{ioanalal04@gmail.com}, \texttt{constantin.morarescu@univ-lorraine.fr}, \texttt{jamal.daafouz@univ-lorraine.fr}, \texttt{lucian@busoniu.net}.}
}

\IEEEaftertitletext{\vskip -0.5em} 

\maketitle
\thispagestyle{empty}
\pagestyle{empty}

\begin{abstract}
We propose two new optimistic planning algorithms for nonlinear hybrid-input systems, in which the input has both a continuous and a discrete component, and the discrete component must respect a dwell-time constraint. Both algorithms select sets of input sequences for refinement at each step, along with a continuous or discrete step to refine (split). The dwell-time constraint means that the discrete splits must keep the discrete mode constant if the required dwell-time is not yet reached. Convergence rate guarantees are provided for both algorithms, which show the dependency between the near-optimality of the sequence returned and the computational budget. The rates depend on a novel complexity measure of the dwell-time constrained problem. We present simulation results for two problems, an adaptive-quantization networked control system and a model for the COVID pandemic. 
\end{abstract}

\section{Introduction}\label{sec:intro}
We consider optimal control of hybrid-input systems in which the discrete input is subject to a minimum dwell-time constraint. A hybrid input has both a continuous and a discrete component, and the dwell-time is the number of steps elapsed before the discrete input changes its value. The dwell-time constraint is motivated by preventing fast switches, either due to physical limitations or to increase performance \cite{heydari2017optimal,allerhand2010robust}. Hybrid-input systems occur e.g.\ in robotics \cite{buss2002nonlinear}, industrial multiple-tanks systems \cite{slupphaug1997mpc} or the automotive industry \cite{lygeros2008hybrid}. Moreover, in networked control systems (NCS), the continuous input can be dynamically quantized \cite{liu2015dynamic}, where the quantization mode is the discrete input. Several methods can be used to solve hybrid-input problems \textit{without} dwell-time constraints, among which branch-and-bound approaches \cite{buss2002nonlinear}, switching control \cite{mareczek1998robust}, or MPC \cite{slupphaug1997mpc}. \textcolor{black}{Optimal control of switched systems is also presented in \cite{zhu2015optimal, kamgarpour2012optimal} (see also references therein), which however do not consider hybrid inputs or dwell-time constraints. For linear dynamics, \cite{DuanWu:14} jointly designs a dwell-time constrained mode sequence and the continuous input.}

Differently from these methods, our focus here is on hybrid-input systems with dwell-time constraints, in which dynamics can be general nonlinear and cost functions arbitrary, as long as both are Lipschitz with respect to the state and the continuous input. The latter input must be scalar, a restriction that can be relaxed at extra computational cost. For such systems, in a first main contribution of this paper (C1), we propose two methods, called OPHIS$\Delta$ and SOPHIS$\Delta$: Optimistic Planning for Hybrid-Input Systems with dwell-time, and Simultaneous OPHIS$\Delta$. Both algorithms produce at each step an open-loop sequence, and are meant to be applied in receding horizon. They are an extension to handle dwell-time constraints of the existing OPHIS and SOPHIS methods \cite{lal22optimistic}. This extension is nontrivial since it impacts the way the computational budget is used and thus also the convergence to a near-optimal solution.


Both algorithms belong to the optimistic planning (OP) class \cite{hren2008optimistic} and iteratively partition the space of infinitely long hybrid-input sequences, by choosing for refinement (splitting) one or several sets at each iteration. Dwell-time constraints in OP were addressed before in problems with only discrete inputs \cite{bucsoniu2017planning}, or autonomous switched systems \cite{heydari2017optimal}. In contrast, here we focus on hybrid-input systems. In the methods proposed, for each chosen set, a time step is also selected, together with the type of split (continuous or discrete). The dwell-time constraint is handled during discrete splits, by checking whether enough steps have passed since the last switch (in which case the discrete input can take any possible value) or the constraint is not yet satisfied (in which case the discrete input must be equal to its previous value). 
In OPHIS$\Delta$, a single set is expanded, one that has the maximum upper bound on the value. SOPHIS$\Delta$ refines any sets that may be optimistic regardless of the Lipschitz constants. Thus, the dependence on these constants is eliminated from set selection, but still remains in step selection. In practice, this gives a performance boost for large horizons. 



The second contribution (C2) is a convergence analysis of (S)OPHIS$\Delta$, driven by a novel complexity measure for the dwell-time constrained problem, which requires analyzing the worst-case complexity. Exploiting this new measure, we tailor results from \cite{lal22optimistic} to find convergence rates of the two methods to the optimal value as computation increases. 

Finally, (C3) simulation results are given for two problems, using SOPHIS$\Delta$. 
First, we consider an NCS framework in which the network can be configured to transmit more or less data. Switching between these modes cannot happen too fast, due to an inability to change the configuration of the network too often. Therefore, a dwell-time constraint is imposed. To exemplify this general NCS framework, we chose an inverted pendulum which must be brought upright. The motor command is the continuous input, while the way in which we quantize this value is the discrete input. Then, we discuss a Susceptible-Infectious-Removed (SIR) model \cite{libotte2020determination}, for pandemic evolution, where SOPHIS$\Delta$ is used to determine an optimal strategy to vaccinate the population and choose the level of quarantine needed. When only the vaccination strategy is given as a discrete control input, we recover the results from \cite{libotte2020determination}, while adding a continuous input to represent the level of quarantine gives better results.

\textcolor{black}{Summarizing, we provide two novel algorithms for optimal control of dwell-time constrained hybrid-input systems, and analyze their relation between computation and near-optimality. The key analytical novelty is the dwell-time constraint on the discrete input, which adds complexity to the structure of the tree expanded compared to \cite{lal22optimistic}, and requires a closer look into how discrete and continuous splits are interspersed. A different complexity measure is therefore obtained than in \cite{lal22optimistic}, leading in turn to different convergence rates. Compared to \cite{bucsoniu2017planning}, the continuous input makes the problem significantly more challenging. Finally, the practical relevance of the algorithms is illustrated in two interesting problems from very different domains.}


Next, Section \ref{sec:prb_stat} formalizes the problem, and Section \ref{sec:algos} describes the two algorithms. The convergence rate analysis is given in Section \ref{sec:analysis}, followed by the simulation results in Section \ref{sec:sim}. Conclusions are presented in Section \ref{sec:conc}.
\section{Preliminaries}\label{sec:prb_stat}
We consider a discrete-time nonlinear hybrid-input system:
\begin{equation}
\textcolor{black}{
    x_{k+1}=f(x_k,u_k), \,\,u_k = [c_k, d_k]^T}
\end{equation}
where $x \in X \subseteq \mathbb{R}^m$ is the state and $u \in U$ is the input, which consists of both a continuous action $c_k \in \mathbb{R}$ and a discrete mode $d_k \in \{0,1,...,p\}, p \in \mathbb{N}$. Thus, $U = \mathbb{R}\times\{0,1,...,p\}$. We define a switch as a change from one value of $d$ to another at consecutive steps. The dwell time constraint $\Delta$ is the number of steps during which the discrete input must remain unchanged after a switch.
We also define a reward function $\rho : X \times U \to \mathbb{R}$, representing immediate performance (negative cost) for each state-action pair $(x_k,u_k)$: $r_{k+1} = \rho(x_k,u_k)$. Given an initial state $x_0$ and an infinitely-long sequence of actions (inputs) $\mathbf{u}_{\infty}=(u_0,u_1,...)$, its infinite-horizon discounted value is:
\begin{equation}\label{discountedreturn}
    v(\mathbf{u}_{\infty}) = \sum_{k=0}^{\infty}\gamma^k\rho(x_k,u_k)
    \vspace{-0.7em}
\end{equation}
with $\gamma \in (0,1)$ the discount factor ($\gamma=1$ is excluded). 

Denote by $\mathbf{S}_\Delta^\infty$ the set of infinitely-long action sequences that respect the dwell-time constraint.
The objective is to find the constrained optimal value $v^{*}_{\Delta}:=\mathrm{sup}_{\mathbf{u}_{\infty} \in \mathbf{S}_\Delta^\infty}v(\mathbf{u}_{\infty})$ and a sequence $\mathbf{u}_{\infty}\in \mathbf{S}_\Delta^\infty$ that achieves $v^{*}_{\Delta}$. Note that generally the constrained optimal value is worse than the unconstrained one, so enforcing a dwell-time constraint may lead to a performance loss. 
We make the following assumptions.
\begin{assumption}\label{ass:ass1}
(i) We have $r_k \in [0,1]$ and $c_k \in [0,1]$.\\
(ii) Both the dynamics and the rewards are Lipschitz with respect to the state and the continuous action, i.e., $\exists L_f, L_p\text{ s.t.\ } \forall x, x^{\prime} \in X \text{ and }c,c^{\prime} \in [0, 1]$:
\begin{equation*}\label{lip}
\begin{aligned}
    \|f(x,[c,d]^T)-f(x^\prime,[c^\prime,d]^T)\| & \leq L_f(\|x-x^\prime\| +|c-c^\prime|)\\
    |\rho(x,[c,d]^T)-\rho(x^\prime,[c^\prime,d]^T)| & \leq L_\rho(\|x-x^\prime\| +|c-c^\prime|)\\
    \end{aligned}
\end{equation*}
(iii) $\gamma L_f < 1$.
\end{assumption}

 Bounded costs like in (i) are typical in e.g.\ reinforcement learning for control \cite{sutton2018reinforcement}, and together with discounting they ensure boundedness of the sequence values. The bounded continuous action is needed because we will numerically refine its interval, \textcolor{black}{and is often naturally satisfied due to physical limitations, while the unit interval can be reached by scaling other intervals. }Note that now $U=([0,1]\times\{0,1,...,p\})$. In (ii), Lipschitz continuity is only imposed w.r.t.~the continuous component $c$ of the action, whereas the variation w.r.t.~$d$ can be arbitrary. Note also that (ii) allows nondifferentiable dynamics and rewards, helping to model e.g. saturations, actuator dead-zones, etc. \textcolor{black}{and is not a greatly restrictive property, since usual dynamics and
cost functions satisfy it.} The relationship in (iii) means that the dynamics should become contractive when combined with a shrink rate equal to $\gamma$. \textcolor{black}{This condition is the strongest among the three; it may be relaxed in future work using stability, see \cite{granzotto2019optimistic}.} 

The next property applies to constrained as well as unconstrained sequences. \vspace{-0.3em}
\begin{lemma}\label{lemma:lemmadrive} \normalfont{\cite{lal22optimistic}} \textit{For any two sequences $\mathbf{u}_{\infty}, \mathbf{u}'_{\infty} \in U^{\infty}$:
\begin{multline}\label{eq:metric}
    |v(\mathbf{u}_{\infty}) - v(\mathbf{u}'_{\infty})| \\
    \leq L_\rho\sum_{k=0}^{D-1}|c_k-c'_k|\gamma^k \frac{1-(\gamma L_f)^{D-k}}{1-\gamma L_f}+\frac{\gamma^D}{1-\gamma}
\end{multline}
where $D$ is the first step $k$ at which $d_k \neq d'_k$.}
\end{lemma}
\vspace{-0.3em}
The two terms on the right-hand side of the inequality correspond to the continuous and discrete actions, respectively.

\section{(S)OPHIS with a dwell-time constraint}\label{sec:algos}
This section introduces two new algorithms for the hybrid-input problem with dwell-time constraints of Section \ref{sec:prb_stat}. These algorithms are generalizations of (S)OPHIS \cite{lal22optimistic}, and simplify to them when the problem is unconstrained ($\Delta=1$). Moreover, when the continuous input does not exist, OPHIS$\Delta$ specializes to OP$\delta$ from \cite{bucsoniu2017planning}. The set and step selection rules, as well as continuous-input refinements, are similar to those for (S)OPHIS. The novelty in (S)OPHIS$\Delta$ is the way in which a discrete split is carried out, which is different depending on whether the minimum dwell-time has been surpassed. This will have non-trivial consequences for the complexity of the algorithms in the analysis.


A set of hybrid inputs consists of a continuous-action interval $\mu$ and a discrete action set $\sigma$ for each step $k$:
\vspace{-0.3em}
\begin{equation}
    \mathbb{S}_i = \prod_{k=0}^\infty(\mu_{i,k},\sigma_{i,k})
    \vspace{-0.4em}
\end{equation}
where $\prod$ means the repeated application of the cross-product $\times$, and notation $(\mu, \sigma)$ means a set in which $c\in \mu$ and $d \in \sigma$. For clarity, from now on we will refer to the set per step $k$, $(\mu_{i,k},\sigma_{i,k})$, as a \textit{pair}, and the infinite-horizon $\mathbb{S}_i$ as a \textit{set}. 

For a set $i$, $D_i$ and $C_i$ are respectively the discrete and continuous horizons (numbers of refined discrete and continuous steps). Any step $k < C_i$ has already been refined and its interval is thus strictly smaller than $[0,1]$, whereas for all $k \geq C_i$, $\mu_{i,k}=[0,1]$. For all $k < D_i$, $\sigma_{i,k} = d_{i,k}$, a single, definite value, and for all $k \geq D_i$, $\sigma_{i,k}=\{0,1,...,p\}$. A sequence of actions in set $i$ is then $(u_{i,0}, u_{i,1}, u_{i,2},...)$, where
$    u_{i,k} = \begin{bmatrix}
       c_{i,k},\,
       d_{i,k}
    \end{bmatrix}^T$ 
and $c_{i,k} \in \mu_{i,k}$, $d_{i,k} \in \sigma_{i,k}$. 

Each set has a corresponding dwell-time $\Delta_i$, equal to the number of steps since the last switch for the discrete input:
\vspace{-0.3em}
\begin{equation}
\Delta_i =  \underset{\Delta^{'}}{\max} \text{ s.t. } d_{i,D_i-\Delta^{''}} = d_{i,D_{i}}, \forall \Delta^{''} \leq \Delta^{'}
    \vspace{-0.3em}
\end{equation}

\textcolor{black}{Consider now reward $r_{i,k+1} = \rho(x_{i,k},u_{i,k})$, where we refer by $c_{i,k}$ to the specific action that is at the center of interval $\mu_{i,k}$. Define then the sample value of a set $i$}:
\vspace{-0.4em}
\begin{equation}
    \textcolor{black}{v(i) = \sum_{k=0}^{D_i-1}\gamma^kr_{i,k+1}}
    \vspace{-0.4em}
\end{equation}
\textcolor{black}{Each continuous interval $\mu_{i,k}$ has a length $a_{i,k}$. For $k \geq C_i$, $a_{i,k}=1$. For each set, we define its diameter $\delta(i)$ in the semimetric of \eqref{eq:metric}, so that}:
\begin{equation}\label{diam}
\begin{array}{cc}
     \textcolor{black}{\underset{\mathbf{u}_{\infty},\mathbf{u}'_{\infty} \in \mathbb{S}_i}{\text{sup}} |v(\mathbf{u}_{\infty}) - v(\mathbf{u}'_{\infty})|  \leq \delta(i)}\\
   \textcolor{black}{\delta(i)= L_{\rho}\sum_{k=0}^{D_i-1}a_{i,k}\gamma^k\frac{1-(\gamma L_f)^{D_i-k}}{1-\gamma L_f}+\frac{\gamma^{D_i}}{1-\gamma}}
\end{array}
\end{equation}
\textcolor{black}{For compactness, denote the contribution of step $k$ in the continuous part of the diameter with $\lambda_k:=L_{\rho}a_k\gamma^k\frac{1-(\gamma L_f)^{D-k}}{1-\gamma L_f}$. }
 
\textcolor{black}{The algorithms work by iteratively (a) selecting sets to refine, (b) choosing a continuous or discrete step to split, and (c) performing a split accordingly. These 3 stages are repeated as long as budget is still available and they are detailed below, followed by an example. For both methods, budget $n$ is the allowed number of calls to $f$ and $\rho$.}

\noindent\textit{\textcolor{black}{(a) Set selection:}}\textcolor{black}{ A set $\mathbb{S}_{i^\dagger}$ is selected for refinement, differently in OPHIS$\Delta$/SOPHIS$\Delta$. 
In OPHIS$\Delta$, given the sample value and diameter of set $i$,  define the upper bound:}
\vspace{-0.3em}
\begin{equation}
    \textcolor{black}{B(i) = v(i) + \delta(i)}
    \vspace{-0.3em}
\end{equation}
\textcolor{black}{so that $v(\mathbf{u}_\infty)\leq B(i), \forall \mathbf{u}_\infty \in \mathbb{S}_i$, which follows from Lemma \ref{lemma:lemmadrive}.}
\textcolor{black}{OPHIS$\Delta$ selects for refinement at each iteration an $optimistic$ set, by maximizing the upper bound:\vspace{-0.4em}}
\begin{equation}\label{eq:setselect}
    \textcolor{black}{i^{\dagger} = \text{arg max}_{i\in \mathbb{A}}B(i)}
\vspace{-0.4em}
\end{equation}
\textcolor{black}{where $\mathbb{A}$ is the collection of all sets created so far.}

\textcolor{black}{In SOPHIS$\Delta$, we eliminate the dependency of the set selection rule on the Lipschitz constants. To this end, we expand at each iteration all sets that may be optimistic for any value of this constant (note however that Assumption \ref{ass:ass1} (ii) is still required). Denote by $H$ the depth in the tree created by both methods, equal to the total number of continuous and discrete expansions done to reach a certain set. Since all sets at depth $H$ have the same shape, their diameters $\delta(i)$ are the same, so the maximum-upper-bound set at that depth can only be a set with the largest value $v(i)$. Thus, at each depth $H$ that still has unexpanded nodes, we expand set $i^\dagger$ with the greatest $v$ value among all sets at that depth. We also configure a maximum depth $H_{\max}(n)$ up to which the expansions are allowed to continue, in order to prevent expanding indefinitely. If $H_{\max}$ grows fast with budget $n$, the algorithm will favour deep searches, whereas a slower growth with $n$ focuses the search on breadth. See the convergence rates in Section \ref{sec:analysis} for more insight on how to select $H_{\max}$.}

\noindent\textit{\textcolor{black}{(b) Step and split type selection:}} \textcolor{black}{After selecting a set $\mathbb{S}_{i^{\dagger}}$, we must choose a step to refine, and decide whether we split discretely or continuously.
To this end, we look at the contribution of each step $k$ up to $D_{i^\dagger}-1$ to the diameter (\ref{diam}), as well as at the contribution $\frac{\gamma^{D_{i^\dagger}}}{1-\gamma}$ of the first unrefined step (at the discrete horizon). Whichever contribution is the greatest dictates where we split. Thus:}
\vspace{-0.4em}
\begin{equation}
    \textcolor{black}{k^{\dagger} = \text{arg max}_{k \in \{0,1,...,D_{i^{\dagger}-1}\}}\lambda_k^{\dagger}}
    \vspace{-0.6em}
\end{equation}
\textcolor{black}{If $\lambda_k^{\dagger} \leq \frac{\gamma^{D_{i^\dagger}}}{1-\gamma}$, we split discretely, at horizon $D_{i^\dagger}$. Otherwise, we make a continuous split, along step $\min(k^{\dagger},C_{i^\dagger})$. }
\textcolor{black}{By this rule, we always have $D_i \geq C_i$ for any set $i$.}

\textcolor{black}{Step selection for SOPHIS$\Delta$ remains the same as for OPHIS$\Delta$, so it will unfortunately still depend on the Lipschitz constant, and there is no way to avoid this. }

\noindent\textit{\textcolor{black}{(c) Performing a split:}} A continuous split can be done along any step $k \leq C_i$, by dividing the interval $\mu_{i,k}$ into $M$ equal pieces and thus generating $M$ new sets.
A discrete split is always done at horizon $D_i$. If $\Delta_i < \Delta$, the dwell-time constraint is active, and only one new set is added, making discrete action $d_k$ definite and equal to $d_{k-1}$, since a switch from the value of $d_{k-1}$ to another is not yet possible. If $\Delta_i \geq \Delta$, the dwell-time constraint is satisfied, any discrete action is eligible, and a discrete split adds $p+1$ new sets that make discrete action $d_k$ definite, one set for each discrete possibility. 
Overall, a tree structure is created, each node being a set. To each set chosen for refinement, children are added corresponding to either distinct intervals for continuous splits, or distinct discrete actions. 

\begin{figure}
    \centering    \includegraphics[width=0.75\columnwidth]{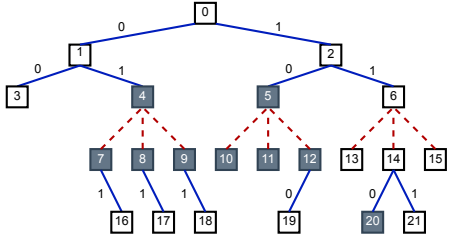}
    \figcapskip
    \caption{Example of tree for $M=3$, $d \in \{0,1\}$, $\Delta=2$}
    \label{fig:tree}
    \figbottomskip
\end{figure}

\noindent\textit{\textcolor{black}{Example:}} An example of the constructed tree is given in Figure \ref{fig:tree}, with $M=3$, $d \in \{0,1\}$, $\Delta=2$. In the figure, the full blue lines correspond to discrete splits, while the red dotted lines mean continuous refinements. For the discrete splits, $0$ and $1$ are added to the branches, to signify the last discrete action. The grey-filled nodes correspond to sets that have $\Delta_i < \Delta$, and are therefore constrained. For those sets, if a discrete split is chosen, the algorithm will only add one child. Any children added by a continuous split of a grey node will also be grey, because continuous action changes do not impact the dwell-time. 

To better understand set splitting, we start by looking at the root node $0$; all continuous intervals are $[0,1]$ and all discrete actions are not yet defined. Then, by a discrete split, we get two new sets, where the first discrete action is defined as $0$ for one set, and $1$ for the other. Sets $1$ and $2$ are now: $\mathbb{S}_1=([0,1],0)\times([0,1],\{0,1\})^\infty$ and $\mathbb{S}_2=([0,1],1) \times ([0,1],\{0,1\})^\infty$. Both sets have $\Delta_i=0$ and are unconstrained, because no switch has yet occurred. Then, set $1$ is split discretely, getting two children nodes $3$ and $4$. Set $3$ is unconstrained, since the last two values of $d$ are equal. However, set $4$ is now constrained, since a switch has been done. Say now that when $4$ is chosen for refinement, a continuous split on step $0$ is done. This adds 3 constrained children, $7$, $8$, and $9$, each corresponding to a third of interval for the first continuous step: $\mathbb{S}_7=([0,1/3],0)\times([0,1],1)\times([0,1],\{0,1\})^\infty$, $\mathbb{S}_8=([1/3,2/3],0)\times([0,1],1)\times([0,1],\{0,1\})^\infty$, $\mathbb{S}_9=([2/3,1],0)\times([0,1],1)\times([0,1],\{0,1\})^\infty$. Since they inherit the discrete input sequence from set $4$, they all have $\Delta_i=1$. Therefore, when a discrete split is done next for sets $7$, $8$, or $9$, only one child corresponding to $d_1=d_0=1$ is added; denote these children by $16$, $17$, and $18$. Now, $\Delta_{16} = \Delta_{17} =\Delta_{18} = 2$, so these nodes are no longer constrained, and the next discrete split from these nodes will again have 2 branches, for $0$ and $1$.

\textcolor{black}{The pseudocode of the two methods is presented in the supplementary material at the end.}
\section{Analysis}\label{sec:analysis}
In this section, we aim to prove that the sub-optimality of the algorithms converges to 0 as budget $n$ increases. Due to the introduction of the minimum dwell-time, the analysis from the unconstrained case \cite{lal22optimistic} cannot be directly applied. We must instead characterize the more complicated tree expanded by the (S)OPHIS$\Delta$ algorithms. To do this, a novel problem complexity measure (branching factor of the constrained tree)  is defined that takes the minimum dwell-time into account. The range of values of this measure is dictated by the largest possible tree. By exploiting the new measure, we are then able to tailor the unconstrained-case results to obtain convergence rates of the new algorithms.


First, recall that we want to find the constrained discounted optimal value $v^*_{\Delta}$, and a dwell-time respecting sequence $i^*$ that approximately achieves it. Similarly to \cite{lal22optimistic}, we have that:
\vspace{-0.9em}
\begin{equation}
    v^*_{\Delta} - v(i^*) \leq \delta_{\min}
    \vspace{-0.4em}
\end{equation}
i.e. the algorithm is near-optimal up to the smallest diameter of any set expanded. Thus, deeper trees give better solutions.

The main part of the analysis focuses on the description of the near-optimal constrained tree, a sub-tree of the full tree that only contains optimistic nodes obeying the dwell-time constraint. This subtree is in fact the one expanded by (S)OPHIS$\Delta$, so characterizing its size is important as it describes the amount of computationrequired.
Define then the set of constrained near-optimal nodes at depth $H$ as:
 \begin{equation}
 \begin{aligned}
     _\Delta\mathcal{T}_H^*=& \{i \text{ at } H \ | \ v(i) + \delta_H\geq v^*_{\Delta}, 
     \\& \text{ and }\forall k < D_i \text{ for which } d_{k-1} \neq d_k,
     \\& d_{k-j}=d_{k-1}, \text{ for } 1 \leq j \leq \Delta, k-j \geq 0\}
\end{aligned}
 \end{equation}
The constrained near-optimal tree has branching factor $\mathcal{K}$:
 \begin{definition}
 The asymptotic branching factor is the smallest $\mathcal{K}$ such that $\exists \mathcal{C} \geq 1$ for which $|_\Delta\mathcal{T}_H^*|\leq\mathcal{C}\mathcal{K}^{H/\Delta}$, $\forall H$, where $|.|$ represents set cardinality.
 \end{definition}
The branching factor is a new measure of complexity for the dwell-time constrained problem. In a simple problem, $\mathcal{K}$ will be small, whereas a complex problem requires expanding many nodes at each depth, leading to a large $\mathcal{K}$. We formalize this next.

\begin{theorem}\label{Kmax}
Branching factor $\mathcal{K}$ is in the range $[\underline{\mathcal{K}}, \overline{\mathcal{K}}]$, where $\underline{\mathcal{K}} = 1$ and  $\overline{\mathcal{K}}=\max(\Delta(p+1), M^\Delta)^{1/ \Delta}$.
\end{theorem}
\vspace*{-1em}\begin{proof}
\noindent The smallest possible value $\underline{\mathcal{K}}=1$ follows easily from the case when the near-optimal tree consists of a single, optimal path.
To find the greatest possible value $\overline{\mathcal{K}}$, we look at the situation in which all rewards are identical, therefore all nodes at each depth have the same $v$ and $B$ and they are all expanded. We want to find the maximum number of nodes at depth $H$, taking into account the $h$ continuous expansions and the discrete number of splits $D$. Recall that a continuous expansion adds $M$ new nodes, and a discrete one adds either $p+1$ new nodes or just $1$, respecting the dwell-time constraint. The number of nodes at depth $H=h+D$ does not depend on the order of continuous and discrete splits. We can therefore consider a different tree, which has first $D$ discrete expansions, all respecting the dwell-time constraint, followed by $h$ continuous splits. This tree will have the same number of nodes at depth $H$ as the original one. Reference \cite{bucsoniu2017planning} proves that
 $
     |_\Delta\mathcal{T}^*_D|\leq (p+1)^2\Delta[\Delta(p+1)]^{D/\Delta}
 $. 
 Therefore:\vspace{-0.2em}
  \begin{equation}
  \begin{array}{rcl}
      |_\Delta\mathcal{T}^*_H| &  \leq & (p+1)^2\Delta[\Delta(p+1)]^{D/\Delta}M^h\\
        &  \leq & (p+1)^2\Delta[\Delta(p+1)]^{H/\Delta}M^H\\
  \end{array}
 \end{equation}
So, $\overline{\mathcal{K}}=\max(\Delta(p+1), M^\Delta)^{1/ \Delta}$, and $\mathcal{K} \in [\underline{\mathcal{K}}, \overline{\mathcal{K}}]$.
\end{proof}
\vspace{-0.5em}
Let us now compare the complexity of the dwell-time-constrained problem with the unconstrained case from \cite{lal22optimistic}. In that setting, the unconstrained tree $\mathcal{T}^*_H$ at depth $H$ has size roughly $m^H$, where the branching factor is $m \in [1, \max(p+1,M)]$. Compared to $m^H$, the constrained tree size $\mathcal{K}^{H/\Delta}$ intuitively emphasizes that discrete choices are made once every $\Delta$ steps. Further, note that the full constrained tree is strictly smaller than the full unconstrained tree, so when expanding full trees using the same budget, the constrained algorithms will reach deeper and have better near-optimality. Now,  this is not immediately visible in the formula, since to extract an easy to interpret branching factor we had to make some conservative replacements (both $h$ and $D$ by $H$). Therefore, to get more insight, consider two cases. In case (i) $\Delta(p+1) \gg M^{\Delta}$, so $|_\Delta\mathcal{T}^*_H| \simeq [\Delta(p+1)]^{H/\Delta}$, significantly smaller than $|\mathcal{T}^*_H| = (p+1)^H$. In other words, since there are many discrete actions, the reduction due to the constraints is significant. Case (ii) $M \gg p+1$, so $|_\Delta\mathcal{T}^*_H| \simeq M^H$, the same as $|\mathcal{T}^*_H|$; since continuous expansions dominate, the reduced number of discrete children is less important. 

The above applies when in both types of problems (constrained and unconstrained), the full tree is expanded. In general, when the two branching factors do not have their maximal values, a clear relationship between the complexity of the constrained and unconstrained problems cannot be established. It could be that the introduction of the constraints makes a constrained-optimal solution easier to distinguish, hence reducing the branching factor/complexity; or, conversely, the constraint could eliminate an optimal solution that would have been easy to find, increasing complexity.


Next, denote $m=\mathcal{K}^{1/\Delta}$, meaning that $m \in [1, \max(\Delta(p+1),M^\Delta)]$. We replace this equivalent branching factor in Theorems 11 and 13 of \cite{lal22optimistic} to get convergence rates of (S)OPHIS$\Delta$, as follows. Recall that $i^*$ denotes the sequence returned by either algorithm; and define $f(n)=\Tilde{O}(g(n))$ to mean that $f(n) \leq a (\log g(n))^b g(n)$ for some $a, b>0$; i.e.\ $f$ behaves like $g$ up to a logarithmic factor. 

\noindent\textit{Convergence rate for OPHIS$\Delta$:} For large budget $n$:\\
a) for $\mathcal{K}>1$: $v^*_{\Delta}-v(i^*)=\Tilde{O}\bigg(\gamma^{\sqrt{\frac{2\overline{\tau}^2(\tau^*-1)\Delta\log n}{\tau^{*^2}\log \mathcal{K}}}}\bigg)$\\
b) for $\mathcal{K}=1$: $v^*_{\Delta}-v(i^*)=\Tilde{O}\Big(\gamma^{n^{1/4}\frac{\overline{\tau}}{\tau^*}\sqrt{\frac{2(\tau^*-1)}{Z\mathcal{C}}}}\Big)$\\
where $\overline{\tau}=\frac{\log(M)}{\log(1/\gamma)}$ and $\tau^*=\lceil{\overline{\tau}}\rceil$.
 

\noindent\textit{Convergence rate for SOPHIS$\Delta$:} For large $n$:\\
a) for $\mathcal{K}>1$, we take $H_{\max}=n^\epsilon$, with $\epsilon \in (0,0.5)$ and we have:\\
$v^*_{\Delta}-v(i^*)=\Tilde{O}\bigg(\gamma\hat\ \Big(\frac{\overline{\tau}}{\tau^*}\sqrt{\frac{(\tau^*-1)(1-2\epsilon)\Delta\log n}{\log \mathcal{K}}}\Big)\bigg)$\\
b) for $\mathcal{K}=1$, we take $H_{\max}=n^{1/3}$, and we have:\\
    $v^*_{\Delta}-v(i^*)=\Tilde{O}\bigg(\gamma \hat\  \Big(n^{1/6}\frac{\overline{\tau}}{\tau^*}\sqrt{2(\tau^*-1)\min \{\frac{1}{\mathcal{C}Z},1\}}\Big)\bigg)$\\
    where $Z=\max(M,p+1)$.

These results say that the sub-optimalities of both OPHIS$\Delta$ and SOPHIS$\Delta$ converge to 0 as $n \to \infty$. The simpler the problem (smaller $\mathcal{K}$), the faster the convergence to 0. In particular, for $\mathcal{K}=1$ convergence is exponential in a power of $n$. For $\mathcal{K}>1$, the multiplication by $\Delta$ at the numerator of the power of $\gamma$ for both algorithms intuitively says that all other things being equal, a larger dwell-time leads to faster convergence. SOPHIS$\Delta$ converges a bit slower than OPHIS$\Delta$, shown by the different powers of $n$ for $\mathcal{K}=1$, and by the appearance of $\epsilon$ for $\mathcal{K}>1$. Further, the results point to a rule for selecting $H_{\max}$: try first with $H_{\max}=n^{1/3}$, and if that does not work well, take $H_{\max}=n^{\epsilon}$ and tune $\epsilon \in (0,0.5)$; for small $\epsilon$, SOPHIS$\Delta$ is nearly as fast as OPHIS$\Delta$. \textcolor{black}{Note that SOPHIS$\Delta$ expands sets for all possible Lipschitz constant values, which intuitively means that it implicitly optimizes the Lipschitz constant for the set selection component. In practice, when larger budgets are available, SOPHIS$\Delta$ is preferred, whereas for smaller budgets, the OPHIS$\Delta$ approach of focusing this limited budget on one value of the Lipschitz constant still pays off.}

\section{Simulation results}\label{sec:sim}
In this section, we present two examples, with the simulations done with SOPHIS$\Delta$, since it provides better results than OPHIS$\Delta$ for long time horizons. The first problem is a quantized NCS framework, applied to an inverted pendulum, and the second a model of the COVID pandemic evolution. \textcolor{black}{For both examples, the algorithm works in receding horizon, so at each step in time, we use it to get an open-loop  sequence of actions, from which we apply the first action.}

\subsection{Quantized NCS framework}
The first problem concerns a Networked Control System (NCS), in which we must transmit commands to an actuator via a network. The precision of the transmitted values is important for performance. By default, many bits are needed in order to transfer a precise value. However, a network with many bits is costly, so we consider sending a more precise value when necessary, and a rougher quantization the rest of the time. This determines a hybrid-input framework, where the value to be transmitted is the continuous input, and the mode of the network (quantization level) is the discrete input, chosen adaptively by the algorithm. Moreover, a dwell time constraint is needed because we cannot switch the configuration of the network too fast.

In the algorithms, we use $M=3$, and therefore, following the tree structure, we use for convenience  trits instead of bits. One trit means the left, center or right interval in a continuous split, and needs 2 bits to be represented. On the actuator's side, a decoding will be made to get the actual control value. If $t$ trits have been sent, this means the interval $[0,1]$ has been split $t$ times.
For example, if we transmit the sequence of $3$ trits: left, center, right, we get first $[0,1/3]$, then $[1/9,2/9]$, and in the end $[5/27,2/9]$. The actual control value will be the center of this last interval, $11/54$.

\begin{figure}
    \centering
    \includegraphics[width=0.86\columnwidth, trim={0.7cm 0.5cm 0.7cm 0},clip]{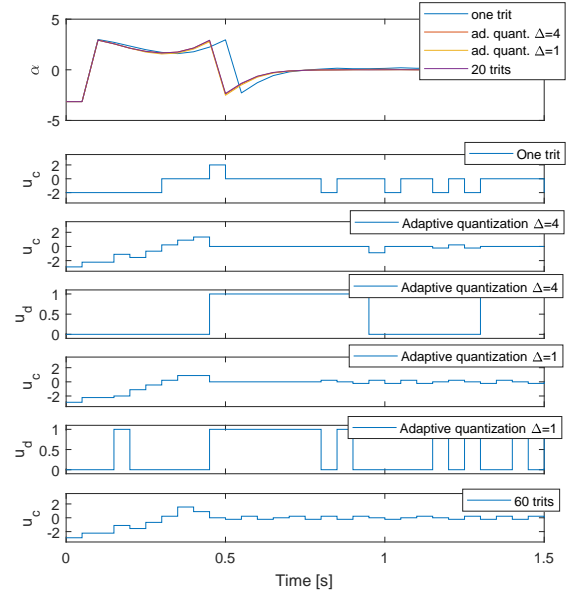}
    \figcapskip
    \caption{Inverted pendulum: State showing the swingup required to reach the target state, quantized continuous inputs, and quantization modes in time}
    \label{fig:states}
\figbottomskip 
\end{figure}

This NCS framework is general: it works for any single-continuous-input system. Next, we apply it to an inverted pendulum, with the nonlinear model given as $\ddot{\alpha}=1/J\cdot[mgl\sin(\alpha)-b\dot{\alpha}-K^2\dot{\alpha}/R+Ku/R]$, with $J=1.91\cdot10^{-4}\mathrm{kgm^2}$, $m=0.055\mathrm{kg}$, $g=9.81\mathrm{m/s^2}$, $l=0.042\mathrm{m}$, $b=3\cdot10^{-6}\mathrm{Nms/rad}$, $K=0.0536\mathrm{Nm/A}$, $R=9.5\mathrm{\Omega}$. We have two states, the angle $\alpha$ and the angular velocity $\dot{\alpha}$. The angle wraps in $[-\pi,\pi]$ and $\dot{\alpha} \in [-15\pi,15\pi]$. The DC motor voltage $u \in [-3,3]V$. The sample time is $T_s = 0.05s$ and we use Euler integration. We use the quadratic reward function $\rho(x_{k+1},u_c) = 1 - 0.75x_{1,k+1}^2/\pi^2-0.25(u_c)^2/9$. We start from $x_0=[-\pi,0]$ (pendulum down) and want to get to $x_f=[0,0]$ (pendulum up). We set our discrete modes to either 0 -- which means sending a sufficiently large number of trits (60) to represent the continuous value after any number of splits made in practice by the algorithms, or 1 -- sending just one trit. We use SOPHIS$\Delta$ with $M=3, L_\rho=1.2, L_f=0.8, \gamma=0.8, n=20000$. As baselines, we look at always sending one trit or 60. Figure \ref{fig:states} shows the states, the applied quantized control voltage, and the modes over time for 4 cases: one trit always, adaptive quantization with dwell-time 4, adaptive quantization with dwell-time 1, and 60 trits always. Constantly using large amount of trits of course gives the best results, but with high network usage. However, the results with adaptive quantization are very similar, and having $\Delta=4$ does not lead to a loss in performance compared to $\Delta=1$, while switching is significantly reduced. Transmitting only one trit all the time degrades the precision of the continuous input and reduces performance. Note that the unquantized input (not shown) is approximately the same as the quantized one, so the algorithm does not waste time refining it more than needed.




\vspace{-0.6em}
\subsection{Susceptible-Infectious-Removed (SIR) model}
\begin{figure}
    \centering
    \includegraphics[width=0.91\columnwidth, trim={0cm 0.6cm 0cm 0.6cm},clip]{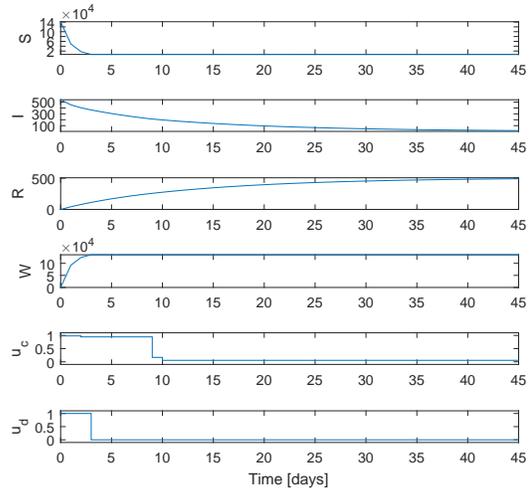}
    \figcapskip
    \caption{SIR: States and inputs in time with dwell-time 2}
    \label{fig:sir_x}
\figbottomskip
\end{figure}

We apply SOPHIS$\Delta$ for a pandemic evolution model, to design the vaccination and quarantine strategy. The model is taken from \cite{libotte2020determination} and its states are: the number of susceptible ($S$), infectious ($I$), removed ($R$) people. In addition, when we include the vaccination control, a new state is introduced: vaccinated ($W$) and the model becomes SIRW. The control variable is $u_d \in \{0,1\}$, with $0$ meaning no vaccines are administered, and $1$ that the maximum percentage of $S$ persons are vaccinated. First, we use OPD \cite{bucsoniu2017planning} (the method OPHIS specializes to when there is no continuous action) just to validate the correctness of our class of methods. As in \cite{libotte2020determination}, we take the values of the parameters: $\beta_{baseline}=0.3566, \gamma=0.0858$ and start from the same initial conditions: $I_{n0}=0.0038, S_{n0}=1-I_{n0}, R_{n0}=0, W_{n0}=0$. With OPD, we recover the same results as \cite{libotte2020determination}, for the mono-objective setting there (achieving the minimum number of infected persons). Then, we add the continuous control variable, equivalent to the level of quarantine. A higher level decreases the infection rate of the virus $\beta$ \cite{pribylova2020seiar}. We set the new value as $\beta_{baseline}-0.5*u_c$. The reward function used is $r=1-0.9998I-0.0001u_c-0.0001u_d$, focusing mostly on reducing the number of infections, but still including small penalties for vaccination (due to its costs) and  quarantine (as it impacts the economy). The simulation results for dwell-time 2 are given in Figure \ref{fig:sir_x}, in which the algorithm works well. Compared to $\int_0^{70}I=8945.42$, the objective function in \cite{libotte2020determination}, we now get $6167.6$. However, recall that we use an additional control variable, which helps reducing the number of infections. In the unconstrained case, we get $5762.5$.\footnote{SOPHIS$\Delta$ also outperforms  \cite{ghosh2019qualitative} on the SIR model in that paper, which natively has 2 control variables. Details are skipped due to space limits.}


\vspace{-0.8em}
\section{Conclusion}
\vspace{-0.6em}\label{sec:conc}
This work presented two new optimistic planning algorithms, OPHIS$\Delta$ and SOPHIS$\Delta$, suited for hybrid-input systems in which the discrete input must respect a dwell-time constraint. The analysis proved that the sub-optimality for each algorithm converges to $0$, as the budget increases. Two simulation examples were given. In the future, we plan to analyze the stability of the algorithms.


\vspace{-1.5em}
\bibliographystyle{style/IEEEtranS}
\bibliography{hop}

\clearpage                                                 
\twocolumn[\begin{center}\large \bf Supplementary material for ``Near-optimal control of nonlinear systems with hybrid inputs and
dwell-time constraints''\vskip 1em \end{center}]

\begin{algorithm}
 \KwIn{state $x_0$, model $f$, reward fcn. $\rho$, split factor $M$, discrete set \{$0,1,...,p$\}, budget $n$, Lipschitz constants $L_f$ and $L_\rho$, discount factor $\gamma$, dwell-time constraint $\Delta$ }
 
 initialize collection of sets $\mathbb{A}$ with $\mathbb{S}_0, D_0= C_0=0$ \\
 \While{budget n not yet exhausted}{
  select set $i^{\dagger} = \text{arg max}_{i\in \mathbb{A}}B(i)$\;
  select dimension with max contribution for continuous actions $k^{\dagger} = \text{arg max}_{k \in \{0,1,...,D_{i^{\dagger}}\}}\lambda_k$\;
  \eIf({\text{  /*split discretely*/}}){$\lambda_{k^{\dagger}} \le \frac{\gamma^{D_{i^\dagger}}}{1-\gamma}$}{
    \eIf(){$\Delta_{i^\dagger}<\Delta$}
    {create one child set from $i^\dagger$\;
    child sets inherit continuous intervals and discrete actions up to dimension $D_{i^\dagger}-1$\;
    action $d_{D_{i^\dagger}}=d_{D_{i^\dagger}-1}$\;
    $\Delta_{child}=\Delta_{i^\dagger}+1$
    }
    {create $p+1$ children sets from $i^\dagger$\;
    create one child set for each $d$ - this action is added for dimension $D_{i^\dagger}$\;
    \eIf(){$d_{i^\dagger}==d_{i^\dagger-1}$}
    {$\Delta_{child}=\Delta_{i^\dagger}+1$\;}{$\Delta_{child}=1$}
    }

   any children will have $D=D_{i^\dagger}+1$ and $C = C_{i^\dagger}$\;
   }({\text{  /*split continuously*/}}){
   expand set $i^\dagger$\ along $k^\dagger$ by creating its $M$ children sets\;
   children sets inherit continuous intervals and discrete actions up to dimension $D_{i^\dagger}-1$\;
   interval at step $k^\dagger$ is refined by splitting into $M$ equal parts\;
   all children will have $D=D_{i^\dagger}$, $\Delta_{child}=\Delta_{i^\dagger}$ and $C = C_{i^\dagger}$ if $k^\dagger \neq  C_{i^\dagger}$, or $C = C_{i^\dagger}+1$ if $k^\dagger =  C_{i^\dagger}$\;
  }
 }
 \KwOut{sequence $\hat{\mathbf{u}}$ of set $i^*=\mathrm{arg}\max_{i \in \mathbb{A}} v(i)$}
 \caption{OPHIS$\Delta$}
\end{algorithm}

\begin{algorithm}
 \KwIn{state $x_0$, model $f$, $\rho$, split factor $M$, discrete set \{$0,1,...,p$\}, budget $n$, Lipschitz constants $L_f$ and $L_\rho$, discount factor $\gamma$, dwell-time constraint $\Delta$, $H_{\max}(n)$}
 
 initialize collection of sets $\mathbb{A}$ with $\mathbb{S}_0, D_0=C_0=0$ \\
 \While{budget still available}{
  $H$ = smallest depth with unexpanded nodes\;
  \eIf(){$H \geq H_{\max}(n)$}{stop and exit the loop;}{
  \While{$H < H_{\max}(n)$}{
  select set $i^{\dagger} = \text{arg max}_{i\in \mathbb{A}}v(i)$\;
  select dimension with max contribution for continuous actions $k^{\dagger} = \text{arg max}_{k \in \{0,1,...,D_{i^{\dagger}}\}}\lambda_k$\;
  \eIf(){$\lambda_{k^{\dagger}} \le \frac{\gamma^{D_{i^\dagger}}}{1-\gamma}$}{
    split discretely (see Alg. 1)
   }(){
   split continuously (see Alg. 1)
  }
  $H=H+1$
  }
 }}
 \KwOut{sequence $\hat{\mathbf{u}}$ of set $i^*=\mathrm{arg}\max_{i \in \mathbb{A}} v(i)$}
 \caption{SOPHIS$\Delta$}
\end{algorithm}

\end{document}